\documentclass[draft]{article}
\usepackage{amsmath,amsfonts,eucal,latexsym}

\ifcase\@ptsize\textheight53\baselineskip\oddsidemargin0.25in\textwidth5.875in
\or\textheight46\baselineskip\oddsidemargin0.125in\textwidth6.125in 
\or\textheight42\baselineskip\oddsidemargin0pt\textwidth6.375in
\fi\evensidemargin\oddsidemargin\marginparwidth0.75in
\newif\iftheoremnumbers\theoremnumberstrue
\newlength{\theoremsep}
\def\@begintheorem#1#2{\par\vskip\theoremsep\penalty\@lowpenalty
      \sl\trivlist\item[\hskip\labelsep{\bf\iftheoremnumbers#2.\ \fi#1.}]%
      \def\condspace{\vskip-\parskip\vspace{0.75ex}}%
      \def\finalcondspace{\vskip-\parskip\vspace{0.75ex}}%
}
\def\@opargbegintheorem#1#2#3{\par\vskip\theoremsep\penalty\@lowpenalty      
      \sl\trivlist\item[\hskip\labelsep{\bf\iftheoremnumbers#2.\ \fi#1\ (#3).}]
      \hskip0em plus.5em minus.4em%
      \def\condspace{\vskip-\parskip\vspace{0.75ex}}%
      \def\finalcondspace{\vskip-\parskip\vspace{0.75ex}}%
      }
\def\@endtheorem{\endtrivlist\vskip-\parskip\pause}
\newlength{\sectsep}
\newlength{\sectitlesep}
\def\section{
\@startsection {section}{1}{\z@}
 {\sectsep}{\sectitlesep}{\large\bf}}
\newlength{\subsectsep}
\newlength{\subsectitlesep}
\def\subsection{
\@startsection{subsection}{2}{\z@}
 {\subsectsep}{\subsectitlesep}{\normalsize}}
\newlength{\titleunit}
\newif\ifmanyauthors
\def\@author#1{\normalsize \lineskip\titleunit #1 \par}
\def\@maketitle{
 \newpage\null\par\vspace*{-6\titleunit} 
 \begin{center}
  \let\@temp@nl\\\let\newline\\
  \def\and{{\footnotesize and}\ }
  \def\@m@acr@{
   \par\footnotesize\vskip\titleunit
   \def\\{\let\and\@@and,\ \def\\{\@temp@nl\def\\{,\ }}}}
  {\Large\bf\@title\par}\vskip3\titleunit
  \let\\\@m@acr@
  \def\@@and{
   \par\normalsize\vskip\titleunit\
             \ifmanyauthors\else{\normalsize and}\par\vskip\titleunit\fi
   \let\\\@m@acr@}
  {\large\@author}
  {\normalsize\it\@date}
 \end{center}
\par \vskip\titleunit
}
\newlength{\myitemsep}\newlength{\mytopsep}\newlength{\mypartopsep}
\let\abstractsize\small
\newlength{\myendabstractskip}
\let\bibsize\small
\def\thebibliography#1{\section*{References\@mkboth
 {References}{References}}\bibsize\list{[\arabic{enumi}]}
 {\itemsep\myitemsep
 \settowidth\labelwidth{[#1]}\leftmargin\labelwidth
 \advance\leftmargin\labelsep
 \usecounter{enumi}}
 \def\newblock{\hskip .11em plus .33em minus .07em}
 \sloppy\clubpenalty4000\widowpenalty4000
 \sfcode`\.=1000\relax}

\def\theequation{\mbox{\thesection;\kern.17em\arabic{equation}}}
\def\@eqnnum{\hbox to.01pt{}\rlap{\rm\hskip-\displaywidth\theequation}}
\leftmargini\parindent
\newlength{\tpu}
\@addtoreset{equation}{section}

\makeatother

\newtheorem{Th}{Theorem}[section]
\newcommand{\pause}{\vskip1.25ex plus.25ex minus.25ex}
\newtheorem{Lemma}[Th]{Lemma}\newtheorem{Cor}[Th]{Corollary}
\newtheorem{Prop}[Th]{Proposition}
\newtheorem{Exm}[Th]{Example}

\newtheorem{Defn}[Th]{Definition}
\newenvironment{Def}{\begin{Defn}\rm}{\end{Defn}}
\newlength{\proofparskip}
\newcommand{\mynewline}[1]{\null\hfill\break\vbox{}
\vskip-\parskip\vskip#1\nopagebreak\noindent}
\newcommand{\condspace}{\vskip-\parskip\vspace{0.5ex}}
\newcommand{\finalcondspace}{\vskip-\parskip}
\newcounter{cond}\let\condformat\arabic
\def\hidecondnumber#1){)}\newcommand{\nocondnumber}{\expandafter\hidecondnumber}
\newenvironment{conditions}[1]{\leavevmode\condspace
\makeconditionlabel{#1}{\condformat{cond}}{cond}{\null}
 \nopagebreak\noindent\begin{list}{\rm\hskip\labelsep\thecond}
 {\leftmargin0pt\labelwidth0pt\usecounter{cond}
 \itemsep\myitemsep\topsep\mytopsep\partopsep\mypartopsep}}
{\finalcondspace\end{list}}
\newenvironment{alphlist}{\begin{conditions}{}
\makeconditionlabel{}{\alph{cond}}{cond}{\null}}{\end{conditions}}

\newenvironment{romanlist}{\makeconditionlabel{}{\roman{cond}}{cond}{\null}
 \leavevmode\condspace\nopagebreak[3]\noindent
 \begin{list}{\makebox[3\labelsep][r]{\thecond}}
 {\leftmargin0pt\labelwidth0pt\usecounter{cond}
 \itemsep\myitemsep\topsep\mytopsep\partopsep\mypartopsep}}
{\finalcondspace\end{list}}
\newcommand{\modifier}{$'$}

\newenvironment{keywords}{\vskip-\parskip\begin{trivlist}\abstractsize
\itemindent1em\item[\textbf{Keywords}]}{\end{trivlist}}
\renewenvironment{abstract}{\begin{trivlist}\abstractsize
\itemindent1em\item[\raisebox{.7ex}{\textbf{Abstract }}]
\mynewline{-3ex}\ignorespaces}{\end{trivlist}}

\newcommand{\makeconditionlabel}[4]{
 \expandafter\def\csname the#3\endcsname{{\rm(#1#2#4)}}}
\newcommand{\starter}{Proof. }
\newcommand{\marker}{ \hspace*{\fill}$\Box$}
\def\restoreqed{\gdef\qed{\marker\global\let\qed\relax}}\restoreqed
\newenvironment{Proof}{\vspace{-2.75ex plus-.5ex minus-.5ex}
 \pagebreak[2]\begin{trivlist}\restoreqed
 \renewcommand{\condspace}{\vspace{-.5ex}}
 \renewcommand{\finalcondspace}{\vskip-2\parskip\vspace{-0.75ex}}
 \setlength{\parskip}{\proofparskip}
 \renewcommand{\newline}{\mynewline{-3ex}}
 \itemindent1em\item[\textsc\starter]}
 {\qed\end{trivlist}\restoreqed
 \vspace{2ex plus.5ex minus.5ex}}

\newcommand{\bbig}[1]{\raisebox{-.1ex}{\large\kern-.1ex$\big#1$\kern-.2ex}}

\renewcommand{\pause}{\vskip.25ex plus1.25ex minus1.25ex}

\author{Anna Avallone \and Paolo Vitolo\\
Dipartimento di Matematica\\ Universit\`a della Basilicata\\
Contrada Macchia Romana\\ 85100 Potenza, Italy
\newline\texttt{avallone@unibas.it}, \texttt{vitolo@unibas.it}}
\date{}
\title{Lattice uniformities on effect algebras}

\newcommand{\n}{}\def\n'{neighbourhood}
\newcommand{\eps}{\varepsilon}
\newcommand{\set}[2]{\{\,#1:#2\,\}}\newcommand{\s}[1]{\{#1\}}
\newcommand{\U}{\mathcal{U}}\newcommand{\V}{\mathcal{V}}
\newcommand{\F}{\mathcal{F}}\newcommand{\G}{\mathcal{G}}
\newcommand{\mathH}{\mathcal{H}}
\newcommand{\card}[1]{\left|#1\right|}
\newcommand{\du}{\mathcal{DU}(L)}\newcommand{\df}{\mathcal{FND}(L)}
\newcommand{\+}{\ifmmode\oplus\else\undefined\fi}
\renewcommand{\-}{\ifmmode\ominus\else\discretionary{-}{}{}\fi}
\let\acuteaccent\'
\renewcommand{\'}{\ifmmode^\perp\else\expandafter\acuteaccent\fi}
\let\circumflex\^
\renewcommand{\^}{\ifmmode\triangle\else\expandafter\circumflex\fi}
\let\hungarianaccent\H
\renewcommand{\H}{\ifmmode\mathH\else\expandafter\hungarianaccent\fi}

\begin{document}
\maketitle
\bibliographystyle{abbrv}

\begin{abstract}
 Let $L$ be a lattice ordered effect algebra. We prove
 that the lattice uniformities on $L$ which make uniformly
 continuous the operations $\-$ and $\+$ of $L$ are
 uniquely determined by their system of neighbourhoods of $0$ and
 form a distributive lattice. Moreover we prove that every such
 uniformity is generated by a family of weakly subadditive
 $[0,+\infty]$-valued functions on $L$. 
\end{abstract}


\section*{Introduction}

Effect algebras have been introduced by D.~J. Foulis
and M.~K. Bennett in 1994 \cite{BF} for modelling unsharp
measurement in a quantum mechanical system. They are a
generalization of many structures which arise in quantum physics
(see~\cite{BC}) and in Mathematical Economics (see~\cite{EZ,BK}), in
particular of orthomodular lattices in non-commutative measure
theory and MV-algebras in fuzzy measure theory. After 1994, there
have been a great number of  papers concerning effect algebras
(see~\cite{DP} for a bibliography).\par
In this paper we study D-uniformities on a lattice ordered effect algebra $L$,
i.e.\ lattice uniformities on $L$ which make uniformly continuous the
operations $\-$ and $\+$ of $L$.\par
The starting point of our paper is observing the key role played by
D-uniformities  in the study of modular measures on $L$ (see~\cite{A01,ABV,AB}),
since every modular measure on $L$ generates a D-uniformity. %, and 
Also of importance is
the role played in the study of modular functions on orthomodular lattices 
(see~\cite{W95}) and of measures on MV-algebras (see~\cite{BW98,Gr00})
by the lattice structure of filters which generate lattice uniformities 
making uniformly continuous the operations of these structures.\par
In the first part of the paper, we give a description of the filters
which are systems of \n's of $0$ in D-uniformities on $L$---called
D-filters---and we prove that there exists an order isomorphism
between the lattice of all D-uniformities on $L$ and
the lattice of all D-filters on $L$. In particular every D-uniformity
is uniquely determined by its system of \n's of $0$.
As a consequence, we obtain that the lattice of all D-uniformities on $L$
is distributive.\par
Our results extend similar results of~\cite{W95} in orthomodular
lattices (see also~\cite{AW})
and of~\cite{BW98} and~\cite{Gr00} in MV-algebras, and give
as particular case the order isomorphism found in~\cite{AV02} between some
lattice congruences and some lattice ideals.\par
In the second part of the paper, we apply the results of the first part
to prove that every D-uniformity on $L$ is generated by a family of
weakly subadditive $[0,+\infty]$-functions on $L$.

\section{Preliminaries}

An \emph{effect algebra}~\cite{DP} is a set $E$, with two distinguished elements 
$0$ and $1$, and a partially defined operation $\+$ such that for all $a,b,c\in E$:
\begin{conditions}{E}
 \item If $a\+b$ is defined, then $b\+a$ is defined and $a\+b=b\+a$.\label{comm}
 \item If $b\+c$ is defined and $a\+(b\+c)$ is defined, then $a\+b$ and 
      $(a\+b)\+c$ are defined, and $a\+(b\+c)=(a\+b)\+c$.\label{assoc}
 \item There exists a unique $a\'\in E$ such that $a\+a\'$ is defined and
      $a\+a\'=1$.\label{ocompl}
 \item If $a\+1$ is defined, then $a=0$.\label{01law}
\end{conditions}
It is easily seen that $a\+0$ is always defined and equals $a$.
If $a\+b$ is defined, we say that $a$ and $b$ are \emph{orthogonal} and write 
$a\perp b$.
% The element $a\'$ given by~\ref{ocompl} is called the 
% \emph{orthocomplement} of $a$.

In an effect algebra $E$ another partially defined operation $\-$
can be defined by the following rule:
\begin{math}%\label{diff}
 c\-a\text{ exists and equals }b\text{ if and only if }
 a\+b\text{ exists and equals }c.
\end{math}
In particular, $a\'=1\-a$. Moreover, if $a\perp b$, then
$a\+b=(a\'\-b)\'=(b\'\-a)\'$.

In an effect algebra $E$ a partial ordering relation $\le$ can be defined
as follows: $a\le c$ if and only if, for some $b\in E$, $a\+b$ exists and
equals $c$. Hence $c\-a$ is defined if and only if $a\le c$.
Moreover $a\perp b$ if and only if $a\le b\'$.

If $a\vee b$ and $a\wedge b$ exist for all $a,b\in E$, then we say that $E$ is
a \emph{lattice ordered effect algebra} (otherwise called \emph{D-lattice}).
In this case, 
we define the \emph{symmetric difference} of any two elements
$a$ and $b$ in $E$ as
\begin{math}
 a\^b=(a\vee b)\-(a\wedge b).%\label{symdiff}
\end{math}
% \pause

Throughout the paper, the symbol $L$ will always denote a
lattice ordered effect algebra. Let us recall that $L$ is an \emph{MV-algebra}
if and only if $(a\vee b)\-b=a\-(a\wedge b)$ for all $a,b\in L$, while
$L$ is an \emph{orthomodular lattice} if and only if $a\'\wedge a=0$
for every $a\in L$.
% \pause

% In the sequel,
We will make use of the following properties (for the proofs
we refer to~\cite{DP}).

\begin{Prop}\label{basic}%\showrefnames
 For all $a,b,c\in L$ we have:
 \begin{romanlist}
  \item\label{D111D2}
       If $a\le b$, then $b\-a\le b$ and $b\-(b\-a)=a$.
  \item\label{D111D3}
       If $a\le b\le c$, then $c\-b\le c\-a$ and $(c\-a)\-(c\-b)=b\-a$.
  \item\label{P112i}
       If $a\le b\le c$, then $b\-a\le c\-a$ and $(c\-a)\-(b\-a)=c\-b$.
%   \item\label{L125}
%        $a\perp b$ if and only if $a\le b\'$, and 
  \item\label{P116ii}
       If $a\le b\'$ and $a\+b\le c$, then $c\-(a\+b)=(c\-a)\-b=(c\-b)\-a$.
  \item\label{P116iii}
       If $a\le b\le c\'$, then $a\+b\le b\+c$ and $(b\+c)\-(a\+c)=b\-a$.
  \item\label{P116iv}
       If $a\le b\le c$, then $a\+(c\-b)=c\-(b\-a)$.
  \item\label{P116viii}
       If $a\le b\'\le c\'$, then $a\+(b\-c)=(a\+b)\-c$.
  \item\label{P18:21,22}
       If $a\le c$ and $b\le c$, then $c\-(a\vee b)=(c\-a)\wedge(c\-b)$
       and $c\-(a\wedge b)=(c\-a)\vee(c\-b)$.
  \item\label{P18:23,4} 
       If $c\le a$ and $c\le b$, then $(a\wedge b)\-c=(a\-c)\wedge(b\-c)$
       and $(a\vee b)\-c=(a\-c)\vee(b\-c)$.
  \item\label{P1862} 
       If $a\le c\'$ and $b\le c\'$, then $(a\vee b)\+c=(a\+c)\vee(b\+c)$
       and $(a\wedge b)\+c=(a\+c)\wedge(b\+c)$.
 \end{romanlist}%\norefnames
\end{Prop}
% \vspace{-2ex}

Let $\U$ be a uniformity on $L$. We say that $\U$ is a \emph{lattice uniformity}
\cite{Web} if the operations $\vee$ and $\wedge$ are uniformly continuous
with respect to $\U$.\par
A \emph{D-uniformity}~\cite{A01} is a lattice uniformity which makes
the operations $\+$ and $\-$ uniformly continuous, too. The set of
all D-uniformities on $L$ will be denoted by $\du$.
It is easy to see that $\du$---ordered by inclusion---is a complete lattice,
with the discrete uniformity and the trivial uniformity as greatest and 
smallest elements, respectively.

Given $U,V\subset L\times L$, we put
\vspace{-1ex}
\begin{align*}
 U\vee V  &=\set{(a_1\vee b_1,a_2\vee b_2)}{(a_1,a_2)\in U,\ (b_1,b_2)\in V},\\
 U\wedge V&=\set{(a_1\wedge b_1,a_2\wedge b_2)}{(a_1,a_2)\in U,\ (b_1,b_2)\in V},\\
 U\-V&=\set{(a_1\-b_1,a_2\-b_2)}
           {b_1\le a_1,\ b_2\le a_2,\ (a_1,a_2)\in U,\ (b_1,b_2)\in V}.
\end{align*}
It is known (see~\cite{Web}) that a uniformity $\U$ on $L$ is a lattice 
uniformity if and only if for every $U\in\U$ there exists $V\in\U$
such that $V\vee\Delta\subset U$ and $V\wedge\Delta\subset U$, where
$\Delta=\set{(a,a)}{a\in L}$.\par
Similarly, it has been shown in~\cite{A01} that a lattice uniformity $\U$ on $L$
is a D-uniformity if and only if for every $U\in\U$ there exists $V\in\U$
such that $V\-\Delta\subset U$ and $\Delta\-V\subset U$.

\section{D-uniformities and D-filters}

\begin{Def}
 A filter $\F$ of subsets of a D-lattice $L$ is called a \emph{D-filter}
 if it satisfies the following:
 \begin{conditions}{F}
  \item$\forall F\in\F\quad\exists F'\in\F:\quad\forall a,b\in F'\quad
        [a\perp b\implies a\+b\in F]$;\label{fsum}
  \item$\forall F\in\F\quad\exists G\in\F:\quad\forall a\in G\quad
        \forall c\in L\quad(a\vee c)\-c\in F$.\label{fweb}
%  \storecounter{cond}{fweb}
 \end{conditions}
 The set of all D-filters on $L$ will be denoted by $\df$.
\end{Def}

Note that, by~\ref{basic}\ref{D111D3}, a filter $\F$ satisfies \ref{fweb}
if and only if, for every $F\in\F$, there exists $G\in\F$ such that,
for all $a\in G$ and all $c\in L$, one has $c\-(a\'\wedge c)\in F$.\par
We shall prove, in Theorem~\ref{filth} below, that $\df$ is isomorphic to 
$\du$ and that $\F$ is a D-filter if and only if $\F$ is the system of
\n's of $0$ in a D-uniformity.

\begin{Lemma}\label{invD}
 For every $a,b,c,d\in L$ such that $c\le a$,
 $c\le b$, $d\ge a$ and $d\ge b$ one has
 $(a\-c)\^(b\-c)=a\^b=(d\-a)\^(d\-b)$.
\end{Lemma}
\begin{Proof}
 Indeed, applying~\ref{basic}\ref{P18:23,4}, %\ref{basic}\ref{P18:23,4}
 and~\ref{basic}\ref{P112i}, one gets
 \begin{math}
  (a\-c)\^(b\-c)
  =\big((a\-c)\vee(b\-c)\big)\-\big((a\-c)\wedge(b\-c)\big)
  =\big((a\vee b)\-c\big)\-\big((a\wedge b)\-c\big)
  =(a\vee b)\-(a\wedge b)=a\^b.
 \end{math}
 Similarly, applying~\ref{basic}\ref{P18:21,22}, %\ref{basic}\ref{P18:21,22},
 and~\ref{basic}\ref{D111D3}, one gets
 \begin{math}
  (d\-a)\^(d\-b)
  =\big((d\-a)\vee(d\-b)\big)\-\big((d\-a)\wedge(d\-b)\big)
  =\big(d\-(a\wedge b)\big)\-\big(d\-(a\vee b)\big)
  =(a\vee b)\-(a\wedge b)=a\^b.
 \end{math}
\end{Proof}

\kern-2ex
\begin{Prop}\label{fprop}%\renewcommand{\finalcondspace}{\vspace*{-3ex}}
 A D-filter $\F$ on $L$ has the following properties:
 \begin{romanlist}
  \item$\forall F\in\F\quad\exists G\in\F:\quad\forall a\in G\quad\forall b\in L 
        \quad[b\le a\implies b\in F]$;\label{0inF}
  \item$\forall F\in\F\quad\exists G\in\F:\quad\forall a,b\in G\qquad
        a\vee b\in F$;\label{unionF}
%   \item$\forall F\in\F\quad\exists G\in\F:\quad\forall a,b\in G\quad
%         a\^b\in F$;\label{DF}
  \item$\forall F\in\F\quad\exists G\in\F:\quad\forall x,y,z\in L\quad 
        [x\^y\in G\implies (x\vee z)\^(y\vee z)\in F]$;
 \label{DvF}
  \item$\forall F\in\F\quad\exists G\in\F:\quad\forall x,y,z\in L\quad 
        [x\^y\in G\implies (x\wedge z)\^(y\wedge z)\in F]$;
 \label{DnF}
  \item$\forall F\in\F\quad\exists G\in\F:\quad\forall x,y,z\in L\quad 
        [x\^y\in G,\ y\^z\in G\implies x\^z\in F]$.
 \label{DDDF}
 \end{romanlist}
\end{Prop}
\begin{Proof}
 \begin{romanlist}
  \item Let $F\in\F$ and let $G\in\F$ such that~\ref{fweb} is satisfied. Given 
       any $a\in G$ and any $b\in L$ with $b\le a$, put $c=a\-b$. Then
       $b=a\-(a\-b)=\big(a\vee(a\-b)\big)\-(a\-b)=(a\vee c)\-c\in F$.
       \pagebreak[2]
  \item Given $F\in\F$, let $F'\in\F$ satisfy~\ref{fsum}, and let $G\in\F$
       satisfy~\ref{fweb} with $F'$ in place of $F$. If $a,b\in G$, then 
       $(a\vee b)\-b\in F'$. Moreover $b\in F'$ by~\ref{0inF}. Therefore 
       $a\vee b=\big((a\vee b)\-b\big)\+b\in F$.\pagebreak[2]
%   \item Given $F\in\F$, let $F_1\in\F$ satisfying~\ref{0inF}, and let $G\in\F$
%        satisfying~\ref{fweb} with $F_1$ in place of $F$. If $a,b\in G$, then 
%        $a\vee b\in F_1$; since $a\^b\le a\vee b$, it follows that
%        $a\^b\in F$.
  \item Let $F\in\F$ and let $G\in\F$ such that~\ref{fweb} is satisfied. Given 
       $x,y,z$ such that $x\^y\in G$, we put $a=x\^y$ and 
       $c=\big((x\vee z)\wedge(y\vee z)\big)\-(x\wedge y)$ and we show that
       $(x\vee z)\^(y\vee z)=(a\vee c)\-c$. First observe that
       $x\vee y\vee z=(x\vee y)\vee\big((x\vee z)\wedge(y\vee z)\big)$.
       Now, applying~\ref{basic}\ref{P1862} and~\ref{basic}\ref{P116iii}, we have:
       \begin{gather*}
         (x\vee z)\^(y\vee z)
         =(x\vee y\vee z)\-\big((x\vee z)\wedge(y\vee z)\big)\\
  =\bbig((x\vee y)\vee\big((x\vee z)\wedge(y\vee z)\big)\bbig)\-
  \big((x\vee z)\wedge(y\vee z)\big)\\
  =\bbig(\big((x\^y)\+(x\wedge y)\big)\vee
  \big((x\vee z)\wedge(y\vee z)\big)\bbig)\-
  \big((x\vee z)\wedge(y\vee z)\big)\\
  =\bbig(\big(a\+(x\wedge y)\big)\vee\big(c\+(x\wedge y)\big)\bbig)\-
  \big(c\+(x\wedge y)\big)\\
  =\big((a\vee c)\+(x\wedge y)\big)\-\big(c\+(x\wedge y)\big)
  =(a\vee c)\-c.
       \end{gather*}
  \item Given $F\in\F$, take $G\in\F$ such that~\ref{DvF} is satisfied, and 
       let $x,y,z$ such that $x\^y\in G$. By Lemma~\ref{invD} we have
       $x\'\^y\'=x\^y$, and therefore
       $(x\wedge z)\^(y\wedge z)=(x\'\vee z\')\'\^(y\'\vee z\')\'=
       (x\'\vee z\')\^(y\'\vee z\')\in F$.
  \item Given $F\in\F$, let $F_1\in\F$ satisfy~\ref{unionF}, let $F_2\in\F$
       satisfy~\ref{0inF} with $F_1$ in place of $F$, let $F_3\in\F$
       satisfy~\ref{fsum} with $F_2$ in place of $F$ and let $G\in\F$
       satisfy~\ref{DvF} with $F_3$ in place of $F$. If $a,b,c\in L$ are such 
       that both $x\^y$ and $x\^z$ belong to $G$, then
       $a=(x\vee y\vee z)\-(y\vee z)=
       \big((x\vee(x\vee z)\big)\^\big(y\vee(y\vee z)\big)\in F_3$
       and $b=(y\vee z)\-z=(y\vee z)\^(z\vee z)\in F_3$ also. It follows 
       that $(x\vee y\vee z)\-z=a\+b\in F_2$, so that $(x\vee z)\-z\in F_1$.
       Similarly one shows that $(x\vee z)\-x\in F_1$. Hence $x\^z=
       \big((x\vee z)\-z\big)\vee\big((x\vee z)\-x\big)\in F$.\qed
 \end{romanlist}
\end{Proof}

\begin{Th}\label{filth}
 \begin{alphlist}
  \item If\/ $\U$ is a D-uniformity, then the filter $\F_\U$
       of \n's of $0$ in $\U$ is a D-filter.\label{unifilt}
  \item Let $\F$ be a D-filter and, for each $F\in\F$, let
       \(F^\^=\set{(a,b)\in L\times L}{a\^b\in F}.\)
       Then $\mathcal B=\set{F^\^}{F\in\F}$ is a base for a D-uniformity
       whose filter of \n's of $0$ is $\F$.\label{filtuni}
  \item The mapping $\Psi\colon\U\mapsto\F_\U$ is an order-isomorphism of
       $\du$ onto $\df$ (both ordered by inclusion).\label{iso}
 \end{alphlist}
\end{Th}
\begin{Proof}
 \begin{alphlist}
  \item Since $\+$ is continuous at $(0,0)$, for every $F\in\F_\U$
       there exists $F'\in\F_\U$ such that if $(a,b)\in F'\times F'$
       and $a\perp b$, then $a\+ b\in F$. This gives~\ref{fsum}.
       To prove~\ref{fweb}, let $F\in\F_\U$ and let $U\in\U$ with
       $U(0)\subseteq F$. By uniform continuity of $\-$ and $\vee$,
       there exist $V_1,V_2\in\U$ such that $V_1\-\Delta\subset U$
       and $V_2\vee\Delta\subset V_1$. Now put $G=V_2(0)$, and consider
       any $a\in G$. Then $(0,a)\in V_2$, so that for every $c\in L$
       we have $(c,a\vee c)\in V_2\vee\Delta\subset V_1$ and hence
       $\big(0,(a\vee c)\-c\big)\in V_1\-\Delta\subset U$ which means
       that $(a\vee c)\-c\in U(0)\subseteq F$.
  \item Clearly $F^\^$ is symmetric and $\Delta\subset F^\^$ for 
       every $F\in\F$. Moreover, given $F_1,F_2\in\F$, let $F_3=\F_1\cap F_2$.
       Then $F_3^\^=F_1^\^\cap F_2^\^$. Finally, if 
       $F\in\F$ and $G\in\F$ satisfies~\ref{fprop}\ref{DDDF}, we have that
       $G^\^\circ G^\^\subseteq F^\^$. Therefore $\mathcal B$
       is a base for a uniformity $\U$.\par
       Now fix $U\in\U$. We show that there exists $V\in\U$ 
       such that both $V\vee\Delta$ and $V\wedge\Delta$ are contained
       in $U$. Let $G\in\F$ satisfy~\ref{fprop}\ref{DvF} and put $V=G^\^$.
       % and~\ref{fprop}\ref{DnF} respectively: we  $V_1=G_1^\^$ 
       Given $(x,y)\in V\vee\Delta$, take $a,b,c\in L$ with $x=a\vee c$,
       $y=b\vee c$ and $(a,b)\in V$, that is $a\^b\in G$. By~\ref{fprop}\ref{DvF},
       we have $x\^y=(a\vee c)\^(b\vee c)\in F$, that is $(x,y)\in F^\^$.
       We conclude that $V_1\vee\Delta\subset F^\^\subseteq U$.
       Since the same $G$ also satisfies~\ref{fprop}\ref{DnF} one sees
       in a similar way that $V\wedge\Delta\subset F^\^\subseteq U$ too.
       Next, we show that there exists $V\in\U$ 
       such that both $V\-\Delta$ and $\Delta\-V$ are contained
       in $U$. Choose $F\in\F$ such that $F^\^\subseteq U$ and put
       $V=F^\^$. By Lemma~\ref{invD}, one has
       \begin{math}
         F^\^\-\Delta=\set{(a\-c,b\-c)}{c\le a,\ c\le b,\ a\^b\in F}
         =\set{(a\-c,b\-c)}{c\le a,\ c\le b,\ (a\-c)\^(b\-c)\in F}=F^\^
       \end{math}
       and similarly one sees that $\Delta\-F^\^=F^\^$.
       Hence $V\-\Delta\subset U$ and $\Delta\-V\subset U$.\par
       It remains to prove that the filter of \n's of $0$ in $\U$ coincides with
       $\F$. First observe that, given any $F\in\F$, we have
       \begin{equation}\label{FD0F}
        F^\^(0)=\set{a\in L}{(0,a)\in F^\^}=\set{a\in L}{a\^0\in F}=F
       \end{equation}
       and therefore $F$ is a \n' of $0$ in $\U$. 
       Conversely, if $G$ is a \n' of $0$ in $\U$, since $\mathcal B$ is a base 
       for $\U$, there exists $F\in\F$ such that $F^\^(0)\subseteq G$.
       By~\eqref{FD0F}, this means that $F\subset G$ and hence $G\in\F$, 
       because $\F$ is a filter.
  \item It follows from~\ref{unifilt} that $\Psi$ maps $\du$ into $\df$. 
       Now for any $\F\in\df$ let $\Phi(\F)$ denote the D-uniformity constructed 
       as in~\ref{filtuni}. Since $\Psi\big(\Phi(\F)\big)=\F$, we have that 
       $\Psi$ is onto. Moreover if $\F_1,\F_2\in\du$ and $\F_1\subset\F_2$, then
       $\set{F^\^}{F\in\F_1}\subseteq\set{F^\^}{F\in\F_1}$ whence
       $\Phi(\F_1)\subseteq\Phi(\F_2)$.
       On the other hand, if $\U_1,\U_2\in\du$ and $\U_1\subset\U_2$, then
       the topology induced by $\U_1$ is coarser than the one induced by $\U_2$,
       hence $\Psi(\U_1)\subseteq\Psi(\U_2)$.\par
       Finally we show that $\Phi=\Psi^{-1}$, so that $\Psi$ is one-to-one.
       Given $\F\in\df$, we consider any $\U\in\du$ such that $\F=\Psi(\U)$
       and prove that $\Phi(\F)=\U$. If $F\in\F$, then it is a \n' of $0$,
       hence there is $U\in\U$ such that $U(0)\subseteq F$. By uniform continuity
       of $\^$, there exists $V\in\U$ with $V\^\Delta\subset U$. Now let
       $(a,b)\in V$. We have $(0,a\^b)=(a\^a,b\^a)\in V\^\Delta\subset U$,
       whence $a\^b\in U(0)\subseteq F$. Hence $V\subset F^\^$ and therefore 
       $\U$ is finer than $\Phi(\F)$. Conversely let $U\in\U$. Consider a
       symmetric $V_1\in\U$ with $V_1\circ V_1\subset U$, and take $V_2,V_3\in\U$
       such that $V_2\vee\Delta\subset V_1$ and $V_3\+\Delta\subset V_2$. 
       Put $F=V_3(0)$, so that $F\in\F$. If $(a,b)\in F^\^$, we have $a\^b\in F$,
       that is $(0,a\^b)\in V_3$. It follows that $(a\wedge b,a\vee b)=
       \big(0\+(a\wedge b),(a\^b)\+(a\wedge b)\big)\in V_3\+\Delta\subset V_2$,
       hence $(a,a\vee b)=\big((a\wedge b)\vee a,(a\vee b)\vee a\big)\in
       V_2\vee\Delta\subset V_1$ and, similarly $(b,a\vee b)\in V_1$. Since
       $V_1^{-1}=V_1$ we also have $(a\vee b,b)\in V$, and then
       $(a,b)\in V_1\circ V_1\subset U$. Therefore $F^\^\subseteq U$.
       We conclude that $\U\subset\Phi(\F)$, whence the equality.\qed
 \end{alphlist}
\end{Proof}\medskip

The reader should note that the above theorem implies, as particular cases,
the results of~\cite[Theor.~2.1]{BW98} and~\cite[Theor.~3.6]{Gr00} for 
MV-algebras, as well as \cite[Theor.~1.1]{W95} for orthomodular lattices.

From Theorem~\ref{filth}\ref{iso}, by restricting to principal filters, one can 
deduce the order isomorphism between D-congruences and D-ideals, which has
been found, using a different approach, in~\cite[Theor.~4.5]{AV02}.

\begin{Prop}
 Let $\F$ be the filter of \n's of $0$ in a D-uniformity $\U$. For every 
 $F\in\F$, let
 \begin{math}
  F^\+=\set{(a,b)\in L\times L}
           {\exists h,k\in F:\ h\perp a,\ k\perp b,\ a\+h=b\-k}
 \end{math}
 and
 \begin{math}
  F^\-=\set{(a,b)\in L\times L}
           {\exists i,j\in F:\ i\le a,\ j\le b,\ a\-i=b\-j}.
 \end{math}
 Then both $\set{F^\+}{F\in\F}$ and $\set{F^\-}{F\in\F}$ are bases for $\U$.
\end{Prop}
\begin{Proof}
 It suffices to show that, for every $F\in\F$, there exist $F_1,F_2\in\F$ 
 such that $F^\+,F^\-\supseteq F_1^\^$ and $F_2^\+,F_2^\-\subseteq F^\^$.\par
 Let $F_1\in\F$ satisfy~\ref{fprop}\ref{0inF}. Given
 $(a,b)\in F_1^\^$, we put $h=(a\vee b)\-a$, $k=(a\vee b)\-b$,
 $i=a\-(a\wedge b)$ and $j=b\-(a\wedge b)$.
 Since $h\le(a\vee b)\-(a\wedge b)=a\^b\in F_1$, we have $h\in F$.
 In the same way one sees that $k$, $i$ and $j$ belong to $F$, too.
 Moreover we have
 \(a\+h=a\+\big((a\vee b)\-a\big)=a\vee b=b\+\big((a\vee b)\-b\big)=b\+k,\)
 so that $(a,b)\in F^\+$.
 Similarly, applying~\ref{basic}\ref{D111D2}, we have
 \(a\-i=a\-\big(a\-(a\wedge b)\big)=a\wedge b=b\-\big(b\-(a\wedge b)\big)=b\-j,\)
 so that $(a,b)\in F^\-$.\par
 Now let $G\in\F$ satisfy~\ref{fprop}\ref{unionF}, and take $F_2\in\F$
 satisfying~\ref{fprop}\ref{0inF} with $G$ in place of $F$. Given
 $(a,b)\in F_2^\+$, there are $h,k\in F_2$ such that $h\perp a$, $k\perp b$ and
 $a\+h=b\+k$. Since $a\vee b\le(a\+h)\vee(b\+k)=a\+h=b\+k$, we get
 $(a\vee b)\-a\le h$ and $(a\vee b)\-b\le k$, so that both
 $(a\vee b)\-a$ and $(a\vee b)\-b$ belong to $G$.
 By~\ref{basic}\ref{P18:21,22}, we have
 $a\^b=\big((a\vee b)\-a\big)\vee\big((a\vee b)\-b\big)$
 hence $a\^b\in F$, i.e.\ $(a,b)\in F^\^$.
 Similarly, given $(a,b)\in F_2^\-$, take $i,j\in F_2$ such that $i\le a$,
 $j\le b$ and $a\-i=b\-j$. Observe that $a\-i=(a\-i)\wedge(b\-j)\le a\wedge b$
 thus, applying~\ref{basic}\ref{D111D2}, $i=a\-(a\-i)\ge a\-(a\wedge b)$.
 It follows that $a\-(a\wedge b)\in G$, and in the same way one sees that
 $b\-(a\wedge b)\in G$, too. By~\ref{basic}\ref{P18:23,4}, we have
 $a\^b=\big(a\-(a\wedge b)\big)\vee\big(b\-(a\wedge b)\big)$
 hence $a\^b\in F$, i.e.\ $(a,b)\in F^\^$.
\end{Proof}

Given $F,G\subset L$, we will put $F\+G=\set{f\+g}{f\perp g,\ f\in F,\ g\in G}$.
Using this notation, condition~\ref{fsum} may be rewritten as follows:
\begin{math}
 \forall F\in\F\quad\exists F'\in\F:\quad F'\+F'\subseteq F.
\end{math}

\begin{Prop}\label{supinf}
 \begin{alphlist}
  \item If $\F,\G\in\df$, then $\set{F\+G}{F\in\F,\ G\in\G}$ is a base for
       % a D-filter $\H$ which is the infimum of $\F$ and $\G$, i.e.\ 
       $\F\wedge\G$ in $\df$.
  \item If $\Gamma\subset\df$, then $\bigvee\Gamma$ in $\df$ is the set of all
       intersections of finite subsets of $\bigcup\Gamma$.
       In particular $\G_1\vee\G_2=\set{G_1\cap G_2}{G_1\in\G_1,\ G_2\in\G_2}$
       for all\/ $\G_1,\G_2\in\df$.
 \end{alphlist}
\end{Prop}
\begin{Proof}
 \begin{alphlist}
  \item First observe that\vadjust{\kern-1.5ex}
       \begin{equation}\label{FpG}
        \forall F\in\F\quad\forall G\in\G\quad F\cup G\subset F\+G.
       \end{equation}
       Indeed, since $0\in G$, one has $F=\set{f\+0}{f\in F}\subseteq
       \set{f\+g}{f\perp g,\ f\in F,\ g\in G}=F\+G$, and similarly for $G$.
       In particular, all sets $F\+G$ with $F\in\F$ and $G\in\G$ are non-empty.
       Now, given $F_1\+G_1$ and $F_2\+G_2$, with $F_1,F_2\in\F$ and 
       $G_1,G_2\in\G$, let $F=F_1\cap F_2$ and $G=G_1\cap G_2$. We have
       \begin{math}
        F\+G=\set{f\+g}{f\perp g,\ f\in F,\ g\in G}\subseteq
 \set{f\+g}{f\perp g,\ f\in F_1,\ g\in G_1}=F_1\+G_1
       \end{math}
       and, similarly, $F\+G\subset F_2\+G_2$. Hence
       $F\+G\subset(F_1\+G_1)\cap(F_2\+G_2)$. Therefore
       $\set{F\+G}{F\in\F,\ G\in\G}$ is a base for a filter which we 
       denote by $\H$.\par
       We prove that $\H$ is a D-filter. Given any $H\in\H$, let $F\in\F$ and
       $G\in\G$ such that $F\+G\subset H$. Take $F',F''\in\F$
       satisfying~\ref{fsum} and~\ref{fweb} respectively, and choose 
       $G',G''\in\G$ in a similar way. Clearly $H'=F'\+G'$ and $H''=F''\+G''$ 
       belong to $\H$. We show that $H'$ satisfies~\ref{fsum} and $H''$ 
       satisfies~\ref{fweb} (with $H$ in place of $F$). If $a$ and $b$ are 
       orthogonal elements of $H'$, then $a=f_1\+g_1$ and $b=f_2\+g_2$,
       where $f_1,f_2\in F'$ and $g_1,g_2\in G'$. Note that $f_1\perp f_2$ and
       $g_1\perp g_2$, hence $f=f_1\+f_2\in F$ and $g=g_1\+g_2\in G$. Therefore
       $a\+b=(f_1\+g_1)\+(f_2\+g_2)=(f_1\+f_2)\+(g_1\+g_2)=f\+g\in F\+G\subset H$.
       Now let $a\in H''$ and $c\in L$. Let $f\in F''$ and $g\in G''$ such that
       $a=f\+g$, and put $d=(f\vee c)\-f$. We have $f'=(f\vee c)\-c\in F$ and
       $g'=(g\vee d)\-d\in G$. Since $g\vee d=g'\+d$ and $f\vee c=f\+d$,
       applying~\ref{basic}\ref{P1862} and~\ref{basic}\ref{P116viii}, we obtain
       \begin{math}
         (a\vee c)\-c=(a\vee f\vee c)\-c=\big((f\+g)\vee(f\vee c)\big)\-c
         =\big((f\+g)\vee(f\+d)\big)\-c=\big(f\+(g\vee d)\big)\-c
         =\big(f\+(g\vee d)\big)\-c=\big(f\+(g'\+d)\big)\-c
         =\big((f\+d)\+g'\big)\-c=\big((f\vee c)\+g'\big)\-c
         =\big((f\vee c)\-c\big)\+g'=f'\+g'\in F\+G\subset H.
       \end{math}\par
       It follows from~\eqref{FpG} that both $\F$ and $\G$ are finer than $\H$. To 
       complete the proof, consider any D-filter $\H'$ such that both $\F$ and $\G$
       are finer than $\H'$. We show that $\H'\subseteq\H$. Let $H\in\H$. 
       By~\ref{fsum}, there exists $H'\in\H'$ such that $H'\+H'\subseteq H$. 
       Since $H'\in\F\cap\G$ we get $H'\+H'\in\H$ and hence $H\in\H$, too.
  \item Let $\F$ be the set of intersections of finite subsets of
       $\bigcup\Gamma$. We show that $\F$ is a filter.\par
       Let $F_1,F_2\in\F$. One has $F_1=\bigcap\F_1$ and $F_1=\bigcap\F_1$, where
       $\F_1$ and $\F_2$ are finite subsets of $\bigcup\Gamma$. If $G=F_1\cap F_2$,
       then $G\in\F$ because it is the intersection of $\F_1\cup\F_2$, which is 
       again a finite subset of $\bigcup\Gamma$. Now let $F\in\F$. Then 
       $F=\bigcap_{i=1}^nF_i$, where $F_i\in\G_i$ and $\G_i\in\Gamma$ for each
       $i\in\s{1,2,\ldots,n}$. If $G\supset F$, let $A=G\setminus F$. 
							For each $i$, one has $G_i=A\cup F_i\in\G_i$, and 
       $\bigcap_{i=1}^nG_i=\bigcap_{i=1}^n(A\cup F_i)=A\cup\bigcap_{i=1}^nF_i=
       A\cup F=G$. Hence $G\in\F$.\par
       Now we check properties~\ref{fsum} and~\ref{fweb}. Let $F\in\F$. As
       above, $F=\bigcap_{i=1}^nF_i$, with $F_i\in\G_i\in\Gamma$. For each $i$, 
       take $F_i'$ and $G_i$ in $\G_i$ satisfying~\ref{fsum} and~\ref{fweb} 
       respectively (with $F_i$ in place of $F$). Put $F'=\bigcap_{i=1}^nF_i'$
       and $G=\bigcap_{i=1}^nG_i$. Clearly $F'$ and $G$ belong to $\F$. We 
       show that $F'$ satisfies~\ref{fsum} and $G$ satisfies~\ref{fweb}.
       If $a$ and $b$ are orthogonal elements of $F'$, then for each
       $i\in\s{1,2,\ldots,n}$ we have $a,b\in F_i'$ and hence $a\+b\in F_i$.
       Therefore $a\+b\in F$. Similarly, if $a\in G$ and $c\in L$, then for each
       $i\in\s{1,2,\ldots,n}$ we have $a\in G_i$ and hence $(a\vee c)\-c\in F_i$.
       Therefore $(a\vee c)\-c\in F$.\par
       Since it is clear that each $\G\in\Gamma$ is contained in $\F$ (indeed 
       every $G$ in $\G$ is the intersection of $\s G$, which a finite subset
       of $\bigcup\Gamma$), it remains to prove that any D-filter which is finer
       than all filters in $\Gamma$ is finer than $\F$, too.
       So let $\G'\in\df$ such that $\G\subset\G'$ for every $\G\in\Gamma$.
       Given $F\in\F$, one has $F=\bigcap_{i=1}^nF_i$ where $F_i\in\G_i\in\Gamma$,
       hence $F_i\in\G'$, for each $i\in\s{1,2,\ldots,n}$. Since $\G'$ is a filter,
       we have $F\in\G'$. We conclude that $\F\subset\G'$.\qed
 \end{alphlist}
\end{Proof}

\begin{Cor}\label{distr}
 $\du$ and $\df$ are distributive (complete) lattices.
\end{Cor}
\begin{Proof}
 By Theorem~\ref{filth}\ref{iso}, it is enough to consider $\df$.
%  Moreover, completeness follows immediately from the previous proposition.
 Let $\F_1$, $\F_2$ and $\G$ be D-filters. We have to verify that
 \begin{math}\label{eqdistr}
  (\F\vee\G_1)\wedge(\F\vee\G_2)\subseteq\F\vee(\G_1\wedge\G_2).
 \end{math}\par
 Given $H\in(\F\vee\G_1)\wedge(\F\vee\G_2)$, take $F_1,F_2\in\F$
 and $G_1,G_2\in\G$ with $(F_1\cap G_1)\+(F_2\cap G_2)\subseteq H$. Put 
 $F=F_1\cap F_2$ and let $F'\in\F$ satisfying~\ref{fprop}\ref{0inF}.
%  We claim that
 We complete the proof by showing that
 \begin{math}
  F'\cap(G_1\+G_2)\subseteq(F_1\cap G_1)\+(F_2\cap G_2).
 \end{math}\par
%  whence~\eqref{eqdistr} follows at once.
 Let $a\in F'\cap(G_1\+G_2)$. Choose $a_1\in G_1$ and $a_2\in G_2$ such 
 that $a=a_1\+a_2$. Since $a_1\le a$ and $a\in F'$, one has $a_1\in F\subset F_1$
 and hence $a_1\in F_1\cap G_1$. Similarly one sees that $a_2\in F_2\cap G_2$.
 Therefore $a=a_1\+a_2\in(F_1\cap G_1)\+(F_2\cap G_2)$.
%  , which proves the claim.
\end{Proof}

\begin{Prop}
 If $\F,\G\in\df$, then $\set{F\wedge G}{F\in\F,\ G\in\G}$ is a base
 for $\F\vee\G$, where $F\wedge G=\set{f\wedge g}{f\in F,\ g\in G}$.
\end{Prop}
\begin{Proof}
 Given $F\in\F$ and $G\in\G$, since
 \(F\cap G=\set{a\wedge a}{a\in F\cap G}
   \subseteq\set{f\wedge g}{f\in F,\ g\in G}=F\wedge G,\)
 it remains to prove that there exist $F'\in\F$ and $G'\in\G$ such that
 $F'\wedge G'\subseteq F\cap G$.
 Take $F'\in\F$ satisfying~\ref{fprop}\ref{0inF}, and let $G'$ be a member of 
 $\G$ satisfying~\ref{fprop}\ref{0inF} also, but with $G$ in place of $F$.
 If $f\in F'$ and $g\in G'$, then $f\wedge g\le f$ hence $f\wedge g\in F$ and,
 similarly, $f\wedge g\le g$ hence $f\wedge g\in G$. Therefore
 $f\wedge g\in F\cap G$.
\end{Proof}

\section{Generating D-uniformities by means of
         \lowercase{$\boldsymbol k$}-submeasures}

\begin{Def}
 Let $k\ge1$. We say that a function $\eta\colon L\to[0,+\infty]$ is
 a \emph{$k$-submeasure} if the following conditions hold:
 \begin{conditions}{S}
  \item$\eta(0)=0$;\label{eta0}
  \item$\forall a,b\in L\quad[a\le b\implies\eta(a)\le\eta(b)]$;\label{monot}
  \item$\forall a,b\in L\quad[a\perp b\implies\eta(a\+b)\le k\eta(a)+\eta(b)]$;
       \label{ksubadd}
  \item$\forall a,b\in L\quad\eta\big((a\vee b)\-b\big)\le k\eta(a)$
       \label{ksubVm}.
 \end{conditions}
 A $1$-submeasure is simply called a \emph{submeasure}.\par
Observe that, if $L$ is an MV-algebra, then every function
$\eta\colon L\to[0,+\infty]$ satisfying~\ref{eta0}, \ref{monot}
and~\ref{ksubadd} with $k=1$ is a submeasure.
\end{Def}

For every $\eps>0$, put $S_\eps=
\set{(x,y)\in\left[0,+\infty\right[\times\left[0,+\infty\right[}
    {\left|x-y\right|<\eps}\cup\s{(+\infty,+\infty)}$.
Then $\set{S_\eps}{\eps>0}$ is base for a uniformity $\mathcal S$ on $[0,+\infty]$
whose relativization to $\left[0,+\infty\right[$ is the usual uniformity, while
$+\infty$ is a uniformly isolated point. In the sequel we will %always
endow $[0,+\infty]$ with this uniformity.

\begin{Prop}\label{submuni}
 For every $k$-submeasure $\eta$ there exists a D-uniformity $\U(\eta)$ which
 is the weakest D-uniformity making $\eta$ uniformly continuous.
\end{Prop}
\begin{Proof}
 For each $\eps>0$, let $F_\eps=\set{a\in L}{\eta(a)<\eps}$. Since
 $F_{\eps_1}\cap F_{\eps_2}=F_{\min\s{\eps_1,\eps_2}}$, the collection
 $\set{F_\eps}{\eps>0}$ is a base for a filter $\F$. We %on $L$. We first
 show that $\F$ is a D-filter. Fix $F$ in $\F$, and take $\eps>0$ with
 $F_\eps\subset F$. Then $F'=F_{\frac\eps{k+1}}$ satisifies~\ref{fsum}
 and $G=F_{\frac\eps k}$ satisfies~\ref{fweb}.\par
 From Theorem~\ref{filth}\ref{filtuni}, the sets $F_\eps^\^$ form a base
 for a D-uniformity $\U(\eta)$. Now we show that $\eta$ is
 $\U(\eta)$-uniformly continuous. Let $\eps>0$ and choose $\delta=\frac\eps k$. 
 For every $(a,b)\in F_\delta^\^$, we have
 \(
  \eta(a\vee b)=\eta\big((a\^b)\+(a\wedge b)\big)\le k\eta(a\^b)+\eta(a\wedge b)
  <\eta(a\wedge b)+k\delta=\eta(a\wedge b)+\eps.
 \)
 Thus, if $\eta(a\vee b)=+\infty$, then $\eta(a\wedge b)=+\infty$
 whence, by monotonicity, $\eta(a)=\eta(b)=+\infty$. Otherwise,
 again by monotonicity, $\eta(a)$ and $\eta(b)$ are both finite, and moreover
 $\left|\eta(a)-\eta(b)\right|\le\left|\eta(a\vee b)-\eta(a\wedge b)\right|<\eps$.
 Hence, in any case, $\big(\eta(a),\eta(b)\big)\in S_\eps$.\par
 Finally, let $\V$ be a D-uniformity on $L$ making $\eta$ uniformly continuous.
 We prove that $\U(\eta)\le\V$, which, by Theorem~\ref{filth}\ref{iso}, is 
 equivalent to $\F\subset\G$, where $\G$ is the filter of \n's of $0$ in $\V$.
 Take any $F\in\F$, and choose $\eps>0$ with $F_\eps\subset F$. Since $\eta$ 
 is continuous at $0$ with respect to $\V$, and $\eta(0)=0$, there is some 
 $G\in\G$ such that if $a\in G$ then $\eta(a)<\eps$, i.e.\ $a\in F_\eps$.
 It follows that $G\subset F_\eps\subset F$, hence $F\in G$.
\end{Proof}

Our aim is to prove a sort of converse of the previous result, namely
Theorem~\ref{genunisub} below.

\begin{Prop}\label{metrmeas}
 Let $k,m\ge1$, and $d$ be a pseudometric %on $L$
 such that for all $a,b,c\in L$:%\par\kern-2.5ex
 \begin{conditions}{P}
  \item$d(a\wedge c,b\wedge c)\le d(a,b)$;\label{kn}
  \item$a\perp c,\ b\perp c\implies d(a\+c,b\+c)\le kd(a,b)$;\label{kplus} 
  \item$d\big((a\vee c)\-c,(b\vee c)\-c\big)\le md(a,b)$;\label{kVm}
  \item$d\big((a\vee c)\-c,0\big)\le kd(a,0)$.\label{kV0}
 \end{conditions}
 For each $a\in L$, put
 \begin{math}
  \tilde\eta(a)=d(a,0).
 \end{math}
 Then $\tilde\eta$ is a $k$-submeasure and $\U(\tilde\eta)$ coincides with the 
 uniformity induced by $d$.
\end{Prop}\par\kern0.5ex
\begin{Proof}
 It is clear that $\tilde\eta$ satisfies \ref{eta0}. Moreover, if $a\le b$, 
 by~\ref{kn} we have
 $\tilde\eta(a)=d(a,0)=d(b\wedge a,0\wedge a)\le d(b,0)=\tilde\eta(b)$
 and~\ref{monot} is proved. Now if $a,b\in L$ are orthogonal,
 then, applying the triangular inequality and~\ref{kplus}, we get
 \(\tilde\eta(a\+b)=d(a\+b,0)\le d(a\+b,b)+d(b,0)\le kd(a,0)+d(b,0)
 =k\tilde\eta(a)+\tilde\eta(b),\)
 that is~\ref{ksubadd}.
 Similarly, taking any $a,b\in L$, by~\ref{kV0} we get
 \(\tilde\eta\big((a\vee b)\-b\big)=d\big((a\vee b)\-b,0\big)
 \le kd(a,0)=k\tilde\eta(a),\)
 that is~\ref{ksubVm}.\par
 Denote by $\V$ the uniformity induced by $d$. The sets
 \(V_\eps=\set{(a,b)\in L\times L}{d(a,b)<\eps}\)
 form a base for $\V$, while the sets
 \(F_\eps^\^=\set{(a,b)\in L\times L}{\tilde\eta(a\^b)<\eps}\)
 form a base for $\U(\tilde\eta)$, as we have seen in Proposition~\ref{submuni}.
 We show that for every $\eps>0$ there exists $\delta>0$ such that
 $F_\delta^\^\subseteq V_\eps$ and $V_\delta\subset F_\eps^\^$.
 This will prove that $\V=\U(\tilde\eta)$.\par
 Take $\delta=\dfrac\eps{2km}$.
 Given $(a,b)\in F_\delta^\^$, applying~\ref{kn}
 and~\ref{kplus}, we have
 \(d(a,b)\le d(a,a\wedge b)+d(a\wedge b,b)
 =d\big((a\vee b)\wedge a,(a\wedge b)\wedge a\big)+
 d\big((a\wedge b)\wedge b,(a\vee b)\wedge b\big)
 \le2d(a\vee b,a\wedge b)=2d\big((a\^b)\+(a\wedge b),0\+(a\wedge b)\big)
 \le2kd(a\^b,0)=2k\tilde\eta(a\^b)<2k\delta\le\eps,\)
 so that $(a,b)\in V_\eps$. Therefore $F_\delta^\^\subseteq V_\eps$.\par
 Now let $(a,b)\in V_\delta$.
 Recall that, by~\ref{basic}\ref{D111D3},
 $(a\^b)\-\big((a\vee b)\-a\big)=a\-(a\wedge b)$
 and, by~\ref{basic}\ref{P18:21,22},
 $\big(a\-(a\wedge b)\big)\wedge\big(b\-(a\wedge b)\big)=0$.
 Hence,
 applying first the triangle inequality and then~\ref{kplus}, \ref{kVm},
 \ref{kn} and again~\ref{kVm},
 we obtain
\begin{math}
  \tilde\eta(a\^b)
  =d(a\^b,0)\le d\big(a\^b,(a\vee b)\-a\big)+d\big((a\vee b)\-a,0\big)%\\
  =d\bbig(\big((a\vee b)\-a\big)\+\big(a\-(a\wedge b)\big),(a\vee b)\-a\bbig)
  +d\big((a\vee b)\-a,(a\vee a)\-a\big)%\\
  \le k d\big(a\-(a\wedge b),0\big)+m d(a,b)%\\
  =k d\bbig(\big(a\-(a\wedge b)\big)\wedge\big(a\-(a\wedge b)\big),
  \big(b\-(a\wedge b)\big)\wedge\big(a\-(a\wedge b)\big)\bbig)+m d(a,b)%\\
  \le k d\big(a\-(a\wedge b),b\-(a\wedge b)\big)+m d(a,b)
  \le km d(a,b)+m d(a,b)<(k+1)m\delta\le\eps
\end{math}
 so that $(a,b)\in F_\eps^\^$. We conclude that 
 $V_\delta\subset F_\eps^\^$.
\end{Proof} 

Recall that if $\mathbb G$ is a topological Abelian group, then a mapping
$\mu\colon L\to\mathbb G$ is called a \emph{modular measure} if the
following hold, for all $a,b\in L$:
\begin{conditions}{M}
 \item$\mu(a)+\mu(b)=\mu(a\vee b)+\mu(a\wedge b)$.\label{modular}
 \item If $a\perp b$, then $\mu(a\+b)=\mu(a)+\mu(b)$.\label{measure}
\end{conditions}%\par
Moreover, (see~\cite[Theor.~3.2]{A01}) the sets
\(\set{(a,b)\in L\times L}{\forall r\le a\^b\quad\mu(r)\in W},\)
where $W$ is a \n' of $0$ in $\mathbb G$, form a base for a D-uniformity $\U$.
This $\U$ is called the D-uniformity \emph{generated by $\mu$}.
% which we denote by $\U(\mu)$.
Note that, in case $\mu$ is positive real-valued (hence in 
particular a submeasure), $\U$ agrees with the $\U(\mu)$ constructed in
Proposition~\ref{submuni}.
\pagebreak[2]

\begin{Th}\label{genunisub}
 Let $\U$ be a D-uniformity on $L$. Then:
 \begin{alphlist}
  \item\label{uniksub} For every $k>1$ there is a family \vadjust{\kern-1.5ex}
       $\s{\tilde\eta_\lambda}_{\lambda\in\Lambda}$ of $k$-submeasures
       with $\U=\sup\limits_{\lambda\in\Lambda}\U(\tilde\eta_\lambda)$.
       Moreover, if\/ $\U$ has a countable base, we can choose $\card\Lambda=1$.
  \item If\/ $\U$ is generated by a modular measure $\mu\colon L\to\mathbb G$,
       where $\mathbb G$ is a topological Abelian group, then there is
       a family $\s{\tilde\eta_\lambda}_{\lambda\in\Lambda}$ of submeasures
       with $\U=\sup\limits_{\lambda\in\Lambda}\U(\tilde\eta_\lambda)$.
  \item\label{mvunisub} If $L$ is an MV-algebra, there is %\vadjust{\kern-1.5ex}
       a family $\s{\tilde\eta_\lambda}_{\lambda\in\Lambda}$ of submeasures
       with $\U=\sup\limits_{\lambda\in\Lambda}\U(\tilde\eta_\lambda)$.
 \end{alphlist}
\end{Th}
\begin{Proof}\nopagebreak
 \begin{alphlist}
  \item For every $a,b\in L$, put $f(a,b)=a\wedge b$, $g(a,b)=(a\wedge b\')\+b$
       and $h(a,b)=(a\vee b)\-b$.
       By~\cite[Prop.~1.1(b)]{W93}, $\U$ has base consisting of sets $U$
       such that, for every $(a,a')\in U$ and every $b\in L$,
       \(\big(f(a,b),\,f(a',b)\big)=\big(f(b,a),\,f(b,a')\big)\in U.\)
       Since $g$ and $h$ are $\U$-uniformly continuous,
       from~\cite[Prop.~1.2]{W93}
       it follows that $\U$ is generated by a family
       $\s{d_\lambda}_{\lambda\in\Lambda}$ of pseudometrics
       (a single pseudometric if $\Lambda$ is countable)
       such that, for every $\lambda\in\Lambda$ and all $a,a',b,b'\in L$:
       \vadjust{\kern-1.5ex}
       \begin{eqnarray*}
          d_\lambda\big(f(a,b),\,f(a',b')\big)&\le&
          d_\lambda(a,a')+d_\lambda(b,b'),\\
          d_\lambda\big(g(a,b),\,g(a',b')\big)&\le&
          k\big(d_\lambda(a,a')+d_\lambda(b,b')\big),\\
          d_\lambda\big(h(a,b),\,h(a',b')\big)&\le&
          k\big(d_\lambda(a,a')+d_\lambda(b,b')\big).
       \end{eqnarray*}
       Clearly each $d_\lambda$ satisfies~\ref{kn} and \ref{kplus},
       as well as \ref{kVm} with $m=k$, hence also \ref{kV0}.
       Therefore, applying Proposition~\ref{metrmeas}, the conclusion follows.
  \item Let $\s{p_\lambda}_{\lambda\in\Lambda}$ be a family of group seminorms 
       generating the topology of $\mathbb G$. By~\cite[Theor.~3]{FT82},
       $\U$ is generated by the family of pseudometrics
       $\s{d_\lambda}_{\lambda\in\Lambda}$ where, for every
       $\lambda\in\Lambda$,
       \[d_\lambda(a,b)
       =\sup\set{p_\lambda\big(\mu(r)-\mu(s)\big)}{r,s\in[a\wedge b,a\vee b]}.\]
       Moreover $d_\lambda$ satisifies~\ref{kn} and the following:
       \begin{equation}\label{kV}
        \forall a,b,c\in L\quad d_\lambda(a\vee c,b\vee c)\le d_\lambda(a,b).
       \end{equation}
       \par
       We can complete the proof, applying Proposition~\ref{metrmeas},
       once we have shown that each $d_\lambda$ satisfies both~\ref{kplus}
       and~\ref{kVm} with $m=k=1$, hence also \ref{kV0}.\par
       Fix $\lambda\in\Lambda$. Given $a,b\in L$, observe first that
       $d_\lambda(a,b)=\sup\set{p_\lambda\big(\mu(r)\big)}{r\le a\^b}$. 
       Now let $c\in L$. By Lemma~\ref{invD} we have
       $\big((a\vee c)\-c\big)\^\big((b\vee c)\-c\big)=(a\vee c)\^(b\vee c)$.
       Therefore, by~\eqref{kV},
       $d_\lambda\big((a\vee c)\-c\big),\big((b\vee c)\-c\big)=
       d_\lambda(a\vee c,b\vee c)\le d_\lambda(a,b)$. %, which gives~\ref{kVm}.
       Finally, if $c\perp a$ and $c\perp b$, then, again by~\ref{invD}, we have
       $(a\+c)\^(b\+c)=a\^b$. Hence $d_\lambda(a\+c,b\+c)=d_\lambda(a,b)$.
 \item Define $f$, $g$ and $h$ as in the proof of~\ref{uniksub}. 
       By~\cite[Prop.~1.5]{W93}, since $g$ is associative and distributive
       with respect to $f$, the uniformity $\U$ has a base consisting of sets $U$
       such that, for every $(a,a')\in U$ and every $b\in L$,
       \(\big(f(a,b),\,f(a',b)\big)=\big(f(b,a),\,f(b,a')\big)\in U\) and
       \(\big(g(a,b),\,g(a',b)\big)=\big(g(b,a),\,g(b,a')\big)\in U.\)
       Moreover $h$ is $\U$-uniformly continuous, and therefore
       from~\cite[Prop.~1.2]{W93}
       it follows that, for any $m>1$, $\U$ is generated by a family
       $\s{d_\lambda}_{\lambda\in\Lambda}$ of pseudometrics
       (a single pseudometric if $\Lambda$ is countable)
       such that, for every $\lambda\in\Lambda$ and all $a,a',b,b'\in L$:
       \begin{eqnarray*}
          d_\lambda\big(f(a,b),\,f(a',b')\big)&\le&
          d_\lambda(a,a')+d_\lambda(b,b'),\\
          d_\lambda\big(g(a,b),\,g(a',b')\big)&\le&
          \big(d_\lambda(a,a')+d_\lambda(b,b')\big),\\
          d_\lambda\big(h(a,b),\,h(a',b')\big)&\le&
          m\big(d_\lambda(a,a')+d_\lambda(b,b')\big).
       \end{eqnarray*}
       Clearly each $d_\lambda$ satisfies~\ref{kn}, \ref{kplus} with $k=1$
       and \ref{kVm}. It remains to show that \ref{kV0} with $k=1$ is 
       satisified, too. Let $a,c\in L$. By~\ref{kn}, we have
       \begin{math}
       d_\lambda\big((a\vee c)\-c,0\big)
       =d_\lambda\big(a\-(a\wedge c),0\big)%\\
       =d_\lambda\bbig(\big(a\wedge\big(a\-(a\wedge c)\big),
                     0\wedge\big(a\-(a\wedge c)\big)\bbig)
       \le d_\lambda(a,0).
       \end{math}\qed
 \end{alphlist}\medskip
\end{Proof}

The reader should note that~\ref{genunisub}\ref{mvunisub} was already proved
in~\cite[Theor.~2.5]{BW98}.

% \bibliography{UniLat}

\end{document}